\documentclass[12pt, twoside,,a4paper draft]{article}
\usepackage[cp1251]{inputenc}
\usepackage{graphicx}
\newcommand{\CopyName}{O.\ P.\ Makarchuk, D.\ M.\ Karvatskyi
 } 
\newcommand{\NAME}{O.\ P.\ Makarchuk, D.\ M.\ Karvatskyi
 } %
\newcommand{\Year}{2024} 
\newcommand{\rightheadtext}{On the Lebesgue measure of one generalised set of subsums of geometric series
} 
\usepackage[english,russian,ukrainian]{babel}
\usepackage{ms_we4}

\renewcommand{\refname}{\refnam}
\vs
\newcommand{\tit}{On the Lebesgue measure of one generalised set of subsums of geometric series} 
\date{05.06.2023}
\begin{document}
\vspace{0.3in}
\markboth{{\NAME}}{{\rightheadtext}}
\begin{center} \textsc {\CopyName} \end{center}
\begin{center} \renewcommand{\baselinestretch}{1.3}\bf {\tit} \end{center}

\vspace{20pt plus 0.5pt} {\abstract{ \noindent O.\ P.\ Makarchuk,  D.\ M.\ Karvatskyi,
\textit{On the Lebesgue measure of one generalised set of subsums of geometric series.}\matref \vspace{3pt} 

In the present paper, we study a set that can be treated as a generalised set of subsums for a geometric series. This object was discovered independently in various mathematical aspects. For instance, it is closely related to various systems of representation of real numbers. The main object of this paper was particularly studied by R. Kenyon, who brought up a question about the Lebesgue measure of the set and conjectured that it is positive. Further, Z. Nitecki confirmed the hypothesis by using nontrivial topological techniques. However, the aforementioned result is quite limited, as this particular case should satisfy a rigid condition of homogeneity. Despite the limited progress, the problem remained understudied in a general framework. 

The study of topological, metric, and fractal properties of the set of subsums for a numerical series is a separate research direction in mathematics. On the other hand, the topic is related to another modern mathematical problem, namely, deepening of the Jessen-Wintner theorem for infinite Bernoulli convolutions and their generalisations. The essence of the problem is to reveal the necessary and sufficient conditions for the probability distribution of a random subsum of a geometric series to be absolutely continuous or singular. 

The Jessen-Wintner theorem guarantees that the distribution is pure (pure discrete, pure singular, or pure absolutely continuous).
Meanwhile, the Levy theorem gives us the necessary and sufficient condition for the distribution to be discrete.
Since the set of subsums for an absolutely convergent series coincides with the set of possible outcomes of the corresponding probability distribution, under certain conditions, it allows us to apply various probability theory techniques for its further investigation. In particular, some techniques help us to prove that the above sets have a positive Lebesgue measure and allow to deepen the Jessen-Wintner theorem under certain conditions.

}} \vsk


\renewcommand{\refname}{\refnam}

\vskip10pt

\section{Introduction}

Let $ \gamma \in (0;1)$, $1<s$ be a fixed positive integer, and $||d_{ij} ||$ be an infinite matrix with $s$ rows and real elements. In addition, let the inequality
\begin{equation}\label{eq1}
|d_{ij}|\le L \; \forall i \in \{0;1;...;s-1\}, j \in N,
\end{equation}
hold for some $L>0$.
In what follows, we consider a set
\begin{equation}
M= \left\{\sum_{n=1}^{+\infty}x_{n}\gamma^{n}\middle| x_{n} \in \{d_{0n}, d_{1n}, ... , d_{(s-1)n}\}, \forall n \in N \right \}.
\end{equation}

Topological, metric, and fractal properties of the set $M$ for some values of $\gamma$ and matrixes $||d_{ij} ||$ were studied in \cite{a8, a9, a1, prkr}.
In particular, \cite{a8} is concerned with case $\gamma=\frac{1}{s}$, $d_{jn}=d_{j}$ for each $j \in \{0;1;...;s-1\}$,  $n\in N$, whereas  $(d_{0}, d_{1}, ... ,d_{s-1})$ is a complete residue system modulo $s$. 
The aforementioned paper also conjectures that $\lambda(M)>0$. This hypothesis was further confirmed in paper \cite{a9}.

We are interested in finding less rigid conditions for matrix $||d_{ij} ||$ under which $M$ has a positive Lebesgue measure.
We introduce a probabilistic approach to solve the above problem.

Let $(\psi _{k})$ be a sequence of independent random variables,  where every $\psi _{k}$ takes values $d_{0k}, d_{1k},..., d_{(s-1)k}$ with probabilities $p_{0k} ,p_{1k},...,p_{(s-1)k}$ respectively. We consider a random variable $\psi$ defined by
$$ \psi =\sum _{k=1}^{\infty }\psi _{k} \gamma^{k}.$$
It is easy to see that $M$ is the set of all possible outcomes of $\psi$. In case the probability distribution of $\psi$ contains an absolutely continuous component (or $\psi$ has a pure absolutely continuous distribution), the set $M$ has a positive Lebesgue measure.

According to the Jessen-Wintner theorem \cite{a3}, $\psi$ has a pure probability distribution. From Levy's theorem \cite{a5}, we get that the distribution of $\psi$ is discrete if and only if
$$W\equiv\prod_{n=1}^{+\infty}\max(p_{0n};p_{1n}; ...;p_{(s-1)n})>0.$$

A large number of articles are devoted to the study of the Lebesgue structure of the distribution of $\psi$, depending on various restrictions for $\psi$ \cite{a2, a6, a7, litv, PMK}. In particular, some necessary and sufficient conditions for $\psi$ to have a singular or absolutely continuous distribution were found in \cite{a8, a6, sal}, provided that $\gamma=\frac{1}{s}$, $d_{jn}=j$ for each $j \in \{0;1;...;s-1\}$, $n \in N$. In the present article, we generalise the results of \cite{a8} and \cite{a9}.

\section{Some condition under which $\lambda(M)>0$.}
\begin{theorem}\label{teo1}
Let $(d_{0n}; d_{1n};...;d_{(s-1)n})$ be a set of integer numbers that forms a complete residue system modulo $s$, for each $n \in N$. If $\gamma=\frac{1}{s}$ and condition \eqref{eq1} holds, then $\lambda(M)>0$.
\end{theorem}
\begin{proof}
Let $p_{jn}=\frac{1}{s}$ for each $j \in \{0;1;...;s-1\}$ and $n \in N$. In this case, the distribution of $\psi$ is continuous ($W=0$). It is clear that two random variables $\psi$ and
$$\psi+\sum_{j=1}^{+\infty}\frac{a_{n}}{s^{n}}=\sum_{j=1}^{+\infty}\frac{\psi_{n}+a_{n}}{s^{n}}$$ 
have the same type of probability distribution for an arbitrary bounded sequence of integer numbers $(a_n)$.
Without loss of generality, let
\begin{equation}\label{eq5}
0=d_{0n}<d_{1n}<...<d_{(s-1)n}\le L,
\end{equation}
hold for every $n \in N$. For an integer parameter $r$ and natural parameter $n$ we define
$$S(r;n)=\sum_{(j_{1};j_{2};...;j_{n})}\prod_{k=1}^{n}p_{j_{k}k},$$
where summation is taken over all possible tuples $(j_{1};j_{2};...;j_{n})$ such that $j_{k} \in \{0; 1;...;s-1\}$ for each $k \in \{1;...;n \}$, and
\begin{equation}\label{eq3}
\frac{r}{s^{n}}=\sum_{k=1}^{n}\frac{d_{j_{k}k}}{s^{k}}.
\end{equation}
If identity \eqref{eq3} does not hold over all possible tuples $(j_{1};j_{2};...;j_{n})$, we put $S(r;n)=0$.

Let's consider a function
$$f_{n}(t)=\prod_{k=1}^{n}\sum_{j=0}^{s-1}p_{jk}t^{d_{jk}s^{-k}}.$$
It is easy to see that
$$f_{n}(t)=\sum_{r}S(r;n)t^{rs^{-n}}.$$
Next, we define a polynomial $h_{n}(t)$, depending on $f_{n}(t)$ as follows
$$h_{n}(t)=f_{n}(t^{s^n})=\sum_{r}S(r;n)t^{r}.$$
By mathematical induction on $n$, we attempt to show that
\begin{equation}\label{eq4}
S(r;n)\le \frac{1}{s^n} \; \forall r \in Z, n \in N. 
\end{equation}
For $n=1$, inequality \eqref{eq4} is obviously true. Further, let \eqref{eq4} hold for $n=k$, i.e.,
$$S(r;k)\le \frac{1}{s^k} \; \forall r \in Z.$$
Let's demonstrate that the inequality is correct for $n=k+1$.
It is clear that
$$h_{k+1}(t)=h_{k}(t^{s^k})\sum_{j=0}^{s-1}p_{j(k+1)}t^{d_{j(k+1)}}, $$
hence, for an assigned $r$, there exists an integer $j$ such that $d_{j(k+1)}\equiv r (mod~s)$ and 
$$S(r;k+1)=S(r;k)p_{j(k+1)} \le \frac{1}{s^k}\cdot \frac{1}{s}=\frac{1}{s^{k+1}}.$$
In what follows, we denote
$$W_{n+1}\equiv\sum_{k=n+1}^{+\infty}\frac{\psi_k}{s^{k}}, $$
and calculate
$$F_{\psi}\left(\frac{r+1}{s^n}\right)-F_{\psi}\left(\frac{r}{s^n}\right)=P\left(\psi \in \left[\frac{r}{s^n};\frac{r+1}{s^n}\right]\right)$$
for positive integers $n$ and $r+1$ such that $\frac{r}{s^n} \in M$.
Taking into account \eqref{eq5}, it is not difficult to see that $\psi \in [\frac{r}{s^n};\frac{r+1}{s^n}]$ if and only if
$$\sum_{k=1}^{n}\frac{\psi_{k}}{s^{k}}=\frac{r-l}{s^{n}},$$
$$W_{n+1}\in \left[\frac{l}{s^n};\frac{l+1}{s^n}\right], $$
where $l \in \{0;1;...; [\frac{L}{s-1}]\}$. By \eqref{eq4}, we obtain
$$F_{\psi}\left(\frac{r+1}{s^n}\right)-F_{\psi}\left(\frac{r}{s^n}\right)=\sum_{l}P\left(\sum_{k=1}^{n}\frac{\psi_{k}}{s^{k}}=\frac{r-l}{s^{n}}\right)P\left(W_{n+1}\in \left[\frac{l}{s^n};\frac{l+1}{s^n}\right]\right)= $$
$$\sum_{l}S(r-l;n)P\left(W_{n+1}\in \left[\frac{l}{s^n};\frac{l+1}{s^n}\right]\right) \le \sum_{l}\frac{1}{s^{n}}P\left(W_{n+1}\in \left[\frac{l}{s^n};\frac{l+1}{s^n}\right]\right)=\frac{1}{s^{n}}.$$
For positive integers $m, k+1$ and $n$ such that $m>k$,  $\frac{m}{s^n} \in M,  \frac{k}{s^n} \in M$ we have
$$F_{\psi}\left(\frac{m}{s^n}\right)-F_{\psi}\left(\frac{k}{s^n}\right)=\sum_{j=1}^{m-k}F_{\psi}\left(\frac{k+j}{s^n}\right)-F_{\psi}\left(\frac{k+j-1}{s^n}\right) \le \frac{m}{s^{n}}-\frac{k}{s^{n}}.$$
Taking into account the above inequality and continuity of $F_{\psi}(x)$, it is easy to check that $F_{\psi}(x)$ satisfies the Lipschitz condition
$$|F_{\psi}(x_{1})-F_{\psi}(x_{2})|\le |x_{1}-x_{2}|\; \forall x_{1}, x_{2} \in R.$$
That means the distribution of $\psi$ is absolutely continuous. According to $P(\psi \in M)=1$, we draw a conclusion that $\lambda(M)>0$.
\end{proof}
\section{The Lebesgue structure of distribution of $\psi$.}
\begin{theorem}\label{teo2}
Let $(d_{0n}; d_{1n};...;d_{(s-1)n})$ be a set of integer numbers that forms a complete residue system modulo $s$, for each $n \in N$. If $\gamma=\frac{1}{s}$ and
\begin{equation}\label{eq6}
Q\equiv\prod_{n=1}^{+\infty}\sum_{j=0}^{s-1}\sqrt{\frac{p_{jn}}{s}}>0,  
\end{equation}
then the distribution of $\psi$ is absolutely continuous.
\end{theorem}
\begin{proof}
Without loss of generality, let condition $\eqref{eq5}$ hold. Let $(\xi _{k})$ be a sequence of independent random variables that take values $d_{0n}, d_{1n},..., d_{(s-1)n}$ with equivalent probability $\frac{1}{s}$. Further, we define a random variable $\xi$ by
$$ \xi =\sum _{k=1}^{\infty }\frac{\xi_{k}}{s^k}.$$

According to theorem $\ref{teo1}$, $\xi$ has an absolutely continuous probability distribution.
We define two sequences of probability spaces $\{ (\Omega _{k} ,A_{k} ,\mu _{k} )\} $ and
$\{ (\Omega _{k} ,A_{k} ,\nu _{k} )\}$ in the following way: let
$\Omega _{k} =\{d_{0k}; d_{1k};...;d_{(s-1)k} \} $, $\sigma $-algebra  $A_{k} $ is defined on the power set of $\Omega _{k} $,
$$\mu _{k} (d_{0k})=p_{0k}, \mu _{k} (d_{1k})=p_{1k}, ... ,\mu _{k} (1)=p_{(s-1)k},$$
$$\nu _{k} (d_{0k})=\frac{1}{s}, \nu _{k} (d_{1k})=\frac{1}{s}, ... ,\nu _{k} (d_{(s-1)k})=\frac{1}{s}.$$
It is clear that $\mu _{k} \ll \nu _{k} $ for any $k\in N$. Let us consider infinite products of probabilistic spaces given by
$$(\Omega ,A,\mu )=\prod\limits _{k=1}^{\infty }(\Omega _{k} ,A_{k} ,\mu _{k}),$$
$$ (\Omega ,A,\nu )=\prod\limits _{k=1}^{\infty }(\Omega _{k} ,A_{k} ,\nu _{k} ). $$
 
Due to Kakutani's theorem \cite{a4},  $\mu \ll \nu $ if and only if
$$ \prod\limits _{k=1}^{\infty }\int\limits _{\Omega _{k} }\sqrt{\frac{d\mu _{k} }{d\nu _{k} } }   d\nu _{k} >0\Leftrightarrow  \prod_{n=1}^{+\infty}\sum_{j=0}^{s-1}\sqrt{\frac{p_{jn}}{s}}>0.$$
According to the theorem statement, the last inequality holds. Let's consider a measurable map $\Omega \stackrel{\varphi }{\longrightarrow}M$
defined as
\[\forall \omega \in \Omega :\varphi (\omega )=\sum\limits _{k=1}^{\infty }\frac{\omega _{k}}{s^k}.\]
We define images $\mu ^{*} $ and $\nu ^{*} $ of measures $\mu $ and $\nu $ under map $\varphi $ in a natural way:
\[\mu ^{*} (E)=\mu (\varphi ^{-1} (E)),\]
\[\nu ^{*} (E)=\nu (\varphi ^{-1} (E)),\]
for an arbitrary Borel subset $E$.

Measure $\mu ^{*} $ coincides with probability measure $P_{\psi} $, while measure $\nu ^{*}$ coincides with probability measure $P_{\xi }$. The last one is equivalent to the Lebesgue measure. From the absolutely continuity of measure $\mu $ with respect to $\nu $ follows the absolutely continuity of measure $\mu ^{*} $ with respect to $\nu ^{*} $. Since $\nu ^{*} \sim \lambda $, then $\eqref{eq6}$ implies absolutely continuity of the distribution of $\psi$.
\end{proof}

\begin{theorem}\label{theor3}
Let matrix $||d_{ij}||$ satisfy condition \eqref{eq1}. If $\gamma=\frac{1}{s}$ and 
the identity $W=0=Q$ holds, then the distribution of $\psi$ is singular.
\end{theorem}
\begin{proof}
Without loss of generality, let condition $\eqref{eq5}$ hold. Further, let 
$$\tau_{k}=\sum_{j=1}^{k}\frac{\psi_{j}}{s^j}, \;\;\eta_{k} = \sum_{n=1}^{+\infty} \frac{\psi_{k+n}}{s^n}.$$ 

It is not difficult to see that $\psi = \tau_{k} + \frac{\eta_{k}}{s^k}$ for any positive integer $k$.
Obviously that the random variable $ \tau_{k} $ is taking values $ c_{1}, c_{2}, ... , c_{s^k} $ with probabilities $ q_{1}, q_{2}, ... , q_{s^k} $ respectively.
Moreover, for an arbitrary $j \in \{1;2;...;s^k \}$ there exists a finite sequence $ (\alpha_{1}, \alpha_{2}, ... , \alpha_{k})$ with $\alpha_i \in \{0;1;...;s-1\}$, $i=\overline{(1, k)}$, such that
$$q_{j}=p_{\alpha_{1}1}p_{\alpha_{2}2} ... p_{\alpha_{k}k}. $$

Next, we have
\[\ F_{\psi} (x)  = P(\tau_{k}+\frac{\eta_{k}}{s^k} < x) = \sum_{j=1}^{s^k} q_{j}\cdot P(c_{j}+\frac{\eta_{k}}{s^k} < x) =\]
\begin{equation}\label{trivial}
 =\sum_{j=1}^{s^k} q_{j}\cdot P(\eta_{k} < s^k \cdot (x-c_{j})) = \sum_{j=1}^{s^k} q_{j} \cdot  F_{\eta_{k}} (s^k \cdot (x-c_{j})).
\end{equation}
 
Assume further that $ F_{\psi} (x)$ is absolutely continuous. Then, for any positive integer $k$, the function $ F_{\eta_{k}} (x) $
must be absolutely continuous as well. Taking into account (\ref{trivial}), we get
\[\  p_{\psi} (x) = \sum_{j=1}^{s^k} q_{j}s^{k} p_{\eta_{k}}(s^k\cdot (x-c{j})) \]
almost everywhere. We will use the following notations:
$$A = \int_{-\infty}^{+\infty} \sqrt{p_{\psi} (x)}dx$$
and
$$ B_{k} = \int_{-\infty}^{+\infty} \sqrt{p_{\eta_{k}}(x)}dx.$$

Since $\sqrt{a+b} \leq \sqrt{a} + \sqrt{b}$ for any non-negative numbers $a$ and $b$, we get
$$\ A = \int_{-\infty}^{+\infty} \sqrt{\sum_{j=1}^{s^k} q_{j}s^k p_{\eta_{k}}(s^k(x-c_{j}))}dx \leq 
\sum_{j=1}^{s^k} \sqrt{q_{j}s^k \int_{-\infty}^{+\infty} p_{\eta_k}\cdot(s^k \cdot (x-c_{j}))dx} =$$
$$\  =\sum_{j=1}^{s^k} \frac{\sqrt{q_{j}s^k}}{s^k}B_{k}  = B_{k}\prod_{n=1}^{k}\sum_{j=0}^{s-1}\sqrt{\frac{p_{jn}}{s}}.$$
One can check that
$$ \sum_{n=1}^{+\infty} \frac{d_{(s-1)n}}{s^n} \le \sum_{n=1}^{+\infty} \frac{L}{s^n} = \frac{L}{s-1}, $$
for any $k \in N$. It implies that $P(\eta_{k} \notin [0;\frac{L}{s-1}])=0$, and thus
$$ \int_{-\infty}^{+\infty} \sqrt{p_{\eta_{k}}(x)}dx = \int_{0}^{\frac{L}{s-1}} \sqrt{p_{\eta_{k}}(x)}dx \leq \int_{0}^{\frac{L}{s-1}} \frac{1+p_{\eta_{k}}(x)}{2}dx = \frac{\frac{L}{s-1}+1}{2}.$$
Further, we observe that
$$\ A \leq \prod_{n=1}^{k}\sum_{j=0}^{s-1}\sqrt{\frac{p_{jn}}{s}} \cdot \frac{\frac{L}{s-1}+1}{2} \to 0 \;\; (k \to \infty), $$
hence $ A = 0$. It follows that $ p_{\psi} (x) = 0$ almost everywhere, contrary to the above assumption.
\end{proof}
Taking into account theorems \ref{teo1} and \ref{teo2}, we can deepen the main result of \cite{a6} as follows:
\begin{corollary}
Let $(d_{0n}; d_{1n};...;d_{(s-1)n})$ be a set of integer numbers that forms a complete residue system modulo $s$, for each $n \in N$. If $\gamma=\frac{1}{s}$, then the random variable $\psi$ has a pure distribution,  moreover:

1) the distribution is discrete if and only if $W>0$;

2) the distribution is singular if and only if $W=0=Q$;

3) the distribution is absolutely continuous if and only if $Q>0$.
\end{corollary}

{\footnotesize

\vsk

Department of Dynamical Systems and Fractal Analysis, Institute of Mathematics of NAS of Ukraine, Tereschenkivska Str. 3, Kyiv, 01601, Ukraine.\\ 
First author email: makolpet@gmail.com\\
Second author email: karvatsky@imath.kiev.ua

\vs
}

\end{document}